\documentclass[12pt]{amsart}

\usepackage[french]{babel}
\usepackage[T1]{fontenc}
\usepackage[latin1]{inputenc}

\usepackage{amsfonts,eufrak}

\usepackage{amssymb, amsmath}

\newcommand{\field}[1]{\mathbb{#1}}
\newcommand{\CC}{\field{C}}
\newcommand{\RR}{\field{R}}

\begin{document}
 
\title[Structure des feuilletages]{Structure des feuilletages k\"ahleriens en courbure semi-n\'egative}

\author{Fr\'ed\'eric TOUZET}



\date{\today}


\selectlanguage{french}

\begin{abstract} Nous \'etudions dans cet article quelques propri\'et\'es des feuilletages (transversalement) k\"ahleriens sur une vari\'et\'e compacte lorsque la forme de Ricci transverse est ``suffisamment'' n\'egative. Nous \'etablissons plus pr\'escis\'ement que l'alg\'ebre de Lie du pseudo-groupe d'holonomie est semi-simple. Ceci fournit un crit\'ere qui assure que les feuilles d'un feuilletage holomorphe \'a classe canonique num\'eriquement triviale sont ferm\'ees. 
 \\
\\
{\sc Abstract.} This paper is concerned with (tranversally) K\"ahler foliations. We proved here that  the holonomy pseudo-group's lie algebra is semi-simple under negativity assumptions of the transverse Ricci tensor. As a consequence, we obtain that the leaves of holomorphic foliations with trivial canonical class are closed submanifolds.

\end{abstract}

\maketitle

\section*{Notations et rappels}

Soit $M$ une vari\'et\'e connexe  munie d'un feuilletage régulier ${\mathcal  F}$ de classe ${\mathcal C}^\infty$. 
\\
D\'esignons respectivement par $TM, T\mathcal F, \nu\mathcal F$ le fibr\'e tangent de $M$, le fibr\'e tangent et normal du feuilletage (identifi\'e \'a un suppl\'ementaire de $T\mathcal F$)  et par $T_{\CC}M, T_{\CC}\mathcal F,\nu_{\CC}\mathcal F$ leurs complexifi\'es respectifs. 

On obtient ainsi une d\'ecomposition du fibr\'e en alg\`ebre ext\'erieure

 $${\Lambda}^i{T_{\CC}^*M}=\bigoplus_{r+s=i}\Lambda^{r,s}{T_{\CC}^*M}$$ o\'u $\Lambda^{r,s}{T_{\CC}^*M}={\Lambda}^r{\nu_{\CC}^*\mathcal F}\otimes \Lambda^s{T_{\CC}^*\mathcal F}$
pour tout entier $i$.\\

Relativement \`a cette bigraduation, on constate que la diff\'erentielle $d$ s'\'ecrit comme somme d' op\'erateurs $d_{1,0},d_{0,1}$ et $d_{2,-1}$ de bidegr\'es respectifs $(1,0),(0,1)$ et $(2,-1)$.\\



Rappelons qu'une forme diff\'erentielle $\omega\in\Omega^i(M):=\Gamma ({\Lambda}^i{T_{\CC}M})$ est dite {\it basique} si elle v\'erifie

   $$i_X\omega=i_Xd\omega=0$$

\noindent pour tout champ de vecteur $X$ tangent au feuilletage.
Cette propri\'et\'e est \'evidemment compatible avec le produit ext\'erieur et la diff\'erentiation, de sorte que l'espace $\Omega^*(M/\mathcal F)$ des formes basiques apparait comme un sous-complexe du complexe de de Rham usuel . Son homologie $H^*(M/\mathcal F)$ est appelée {\it cohomologie basique} de $\mathcal F$. 
\\

Par la suite, on s'interessera aux cas des feuilletages {\it transversalement holomorphe}. Dans cette situation, $\mathcal F$ est défini par un cocycle feuillet\'e $\{U_i,f_i,\gamma_{ij}\}$, o\`u $(U_i)$ est un recouvrement ouvert de $M$, $f_i:U_i\rightarrow S$ une submersion au dessus d'une vari\'et\'e holomorphe transverse $S$ de dimension complexe $n$ et $\gamma_{ij}$ un biholomorphisme local tel que sur $U_i\cap U_j$, on ait $f_i=\gamma_{ij}\circ f_j$.\\
Le fibr\'e $\nu\mathcal F$ admet une structure complexe naturelle h\'erit\'e de celle d\'efinie sur l'espace tangent de $M$. Cette structure induit, avec les notations usuelles, une d\'ecomposition en sous-espaces propres sur le complexifi\'e de $\nu^*\mathcal F$,

$$ {\nu_{\CC}^*\mathcal F}={N^*\mathcal F}^{1,0}\oplus{N^*\mathcal F}^{0,1}$$

et d\'etermine ainsi  un scindage du fibr\'e ${\Lambda}^r\nu_{\CC}^*\mathcal F$ sous la forme

  $${\Lambda}^r\nu_{\CC}^*\mathcal F=\bigoplus_{p+q=r}{{N^*\mathcal F}}^{p,q}$$

o\`u ${{N^*\mathcal F}}^{p,q}=\Lambda^p{N^*\mathcal F}^{1,0}\otimes\Lambda^q{N^*\mathcal F}^{0,1}$.\\

  Comme dans le cas ordinaire,   $d_{1,0}$ se d\'ecompose comme somme de deux op\'erateurs $d_{1,0}=\partial  +\overline {\partial }$ tels que 
   
            $$\partial\Gamma({N^*\mathcal F}^{p,q}\otimes\Lambda^s{T_{\CC}^*\mathcal F})\subset \Gamma({{N\mathcal F}^*}^{p+1,q}\otimes\Lambda^s{T_{\CC}^*\mathcal F})$$ et $$\overline {\partial }\Gamma({{N^*\mathcal F}}^{p,q}\otimes\Lambda^s{T_{\CC}^*\mathcal F})\subset \Gamma({{N^*\mathcal F}}^{p,q+1}\otimes\Lambda^s{T_{\CC}^*\mathcal F}).$$


Dans l'espace des formes basiques, la bigraduation pr\'ec\'edente devient

       $$\Omega^r(M/\mathcal F)= \bigoplus_{p+q=r}  \Omega^{p,q}(M/\mathcal F)$$

avec $\Omega^{p,q}(M/\mathcal F)=\Omega^r(M/\mathcal F)\cap \Gamma({{N^*\mathcal F}}^{p,q})$.\\

Soit $\mathcal F$ transversalement holomorphe; on dira que $\mathcal F$ est {\it hermitien} s'il existe sur $\nu\mathcal F$ une m\'etrique hermitienne $g$ invariante par holonomie.
Soit $n$ la codimension complexe de $\mathcal  F$. Il est \'etabli dans (voir \cite{eh})  qu'on a, lorsque la vari\'t\'e ambiante $M$ est {\it compacte}, l'alternative suivante:\\

a) $H^{2n}(M/\mathcal F)=\{0\}$\\

ou\\

b)  $H^{2n}(M/\mathcal F)=\CC$.\\

Un feuilletage admettant la propri\'et\'e b) est dit {\it homologiquement orientable}. Dans ce cas, la classe (en cohomologie {\it basique})de la
forme volume transverse induite par $g$  est un g\'en\'erateur de $H^{2n}(M/\mathcal F)$\\ 

Soit $(\mathcal F,g)$ un feuilletage hermitien. Dans un syst\`eme  de coordonn\'ees holomorphes $z=(z_1,...,z_n)$ qui param\'etre l'espace local des feuilles, la m\'etrique transverse s'\'ecrit
$$g= \sum_{i,j}g_{ij}(z)dz_i d\overline{z_j}$$
avec les conditions $g_{ij}=\overline{g_{ji}}$;
Par analogie au cas usuel, on peut lui  associer sa forme fondamentale $\eta $ et sa forme (o\'u courbure) de Ricci $\rho$. On a $\eta,\rho\in \Omega^{1,1}(M/\mathcal F)$ et localement:

$$\eta=\sqrt{-1}/2\sum_{i,j}g_{ij}dz_i\wedge d\overline{z_j}$$

$$\rho=-\sqrt{-1}/\pi\partial\overline{\partial} Log(\eta^n/{|dz_1\wedge ...\wedge dz_n|}^2)$$

Nous serons par la suite amen\'e \`a consid\'erer l'op\'erateur 

  $$\overline{*}:\Omega^*(M/\mathcal F)\rightarrow \Omega^{2n-*}(M/\mathcal F)$$ attach\'e \'a la m\'etrique $g$ ainsi que les laplaciens basiques associ\'es, $\delta=dd^{\overline{*}}+d^{\overline{*}}d$ et similairement $\Delta_{\partial}$, $\Delta_{\overline{\partial}}$ (pour une d\'efinition pr\'ecise, le lecteur pourra par exemple consulter \cite{elk} ou \cite{eg}).\\
          
  On focalisera notre attention sur les feuilletages {\it k\"ahleriens}, c'est- \'a dire les feuilletages hermitiens tels que la forme fondamentale $\eta$ est {\it ferm\'ee}. Dans cette situation, on a, comme dans le cas classique les relations $\Delta=2\Delta_{\partial}=2\Delta_{\overline{\partial}}$.\\



 D\'esignons par $\mathfrak S$ des germes de transformations infinit\'esimales de $\mathcal F$, par le sous-faisceau  $\mathfrak T$ des germes champs de vecteurs tangents \`a $\mathcal F$ et par $\Xi:=\mathfrak S/\mathfrak T$ le faisceau des germes de champs {\it basiques}.\\
Soit $\mathcal F$ un feuilletage hermitien (ou plus g\'en\'eralement riemannien) sur  M compacte. On note $\mathcal C$ son {\it faisceau transverse central}. Rappelons qu'il s'agit d'un sous-faisceau localement constant de  $\Xi$ d'alg\`ebres de Lie de dimension finie, introduit par Molino (\cite{mo})  et dont les sections locales sont des  champs de Killing  pour la m\'etrique transverse $g$.

Les orbites de ces champs d\'ecrivent l'adh\'erence des feuilles de ${\mathcal F}$ au sens suivant:  chaque point $m\in M$ admet un voisinage ouvert $V$ tel que la restriction  $\mathcal C_{|V}$ soit un faisceau {\it constant} et tel que l'orbite de $m$ (plus exactement sa projection sur l'espace des feuilles local $V/\mathcal F$)    soit ouverte dans(la projection de) l'adh\'erence de la feuille passant par $m$.

Ces germes de champs de Killing peuvent \^etre caract\'eris\'es comme suit: les \'el\'ements du pseudo-groupe d'holonomie suffisamment proches de l'identit\'e sont ceux de la forme $\exp\ \xi$ o\`u $\xi$ est une section locale de $\mathcal C$;
ce faisceau d'alg\`ebre de Lie conf\`ere ainsi au pseudo-groupe d'holonomie une structure de {\it pseudo-goupe de Lie} (\cite{mo}, appendice d'E.Salem).\\

En particulier, lorsqu'on a affaire \`a un feuilletage hermitien, ces champs sont des parties r\'eelles de champs holomorphes. En d'autres termes, si $X$ d\'esigne un tel champ de Killing, on peut \'ecrire $X-iJX$ sous la forme

    $$\sum a_l(z){\partial\over \partial z_l}$$ dans un syst\'eme de coordonn\'es transverses holomorphes $z=(z_1,...,z_n)$, les $a_l$ \'etant holomorphes.\\

On peut facilement identifier le faisceau $\mathcal C$ sur l'exemple suivant.\\
Soit $(M_i,g_i)\ i=1,2$ deux vari\'et\'es hermitiennes compl\`etes. Consid\'erons alors la vari\'et\'e hermitienne $M_1\times M_2$ munie de la m\'etrique produit et supposons qu'il existe un sous-goupe $\Gamma$ d'isom\'etries {\it holomorphes} dont l'action est diagonale, libre, discr\`ete et cocompacte. La vari\'et\'e quotient $M_1\times M_2/\Gamma$ est alors munie de deux feuilletages holomorphes en position transverse correspondant aux facteurs $M_i\ i=1,2$ (respectivement horizontal et vertical). Prenons par exemple le cas du feuilletage horizontal ${\mathcal F}$; l'action de $\Gamma$ \'etant diagonale, on peut \'etudier sa projection $\Gamma_2$ sur le second facteur $M_2$.  Remarquons que $\Gamma_2$ n'est pas n\'ecessairement ferm\'e dans le groupe des isom\'etries $Isom(M_2)$ de $M_2$,  en revanche, son adh\'erence est un groupe de Lie dont l'alg\`ebre de Lie $\mathcal C_0$ est form\'ee de champs de Killing (pour la m\'etrique $g_2$). Le faisceau transverse central associ\'e \'a $\mathcal F$ n'est alors rien d'autre que $\mathcal C=\pi_*(\mathcal C_0)$ o\'u $\pi:M_1\times M_2\rightarrow M_1\times M_2/\Gamma$ d\'esigne l'application de rev\^etement.\\

Concluons cette section par une derni\'ere d\'efinition.\\

Soit $\mathcal F$ un feuilletage transversalement holomorphe.
 Une $(1,1)$ forme r\'eelle basique $\xi$ sera dite {\it semi- n\'egative} si pour toute section $v$ de ${N\mathcal F}^{1,0}$, $\sqrt{-1}\xi(v\wedge \overline{v})\leq 0$ et quasi-n\'egative s'il existe de plus $v$ tel que $\sqrt{-1}\xi(v\wedge \overline{v})$ soit strictement n\'egative en au moins un point de $M$. Nous retrouvons bien s\^ur la d\'efinition habituelle  lorsque $\mathcal F$ est de dimension nulle!\\

\bigskip 
\section{ Pr\'esentation des r\'esultats et commentaires}
\vskip 20 pt
Sauf mention du contraire, $(\mathcal F,g)$ d\'esignera dor\'enavant un feuilletage {\it k\"ahlerien} sur une vari\'et\'e $M$ compacte; $g$ est la m\'etrique kahlerienne transverse et on note respectivement $\eta$ et $\rho$ sa forme fondamentale et sa forme de Ricci.
\medskip

{\bf  Theor\'eme 1.1}-{\it Supposons que  $({\mathcal F}, g)$  est homologiquement orientable et que de plus  la  courbure de Ricci $\rho$ est cohomologue (dans  $\Omega^2(M/\mathcal F)$) \`a une $(1,1)$ forme  basique quasi-n\'egative; alors le faisceau transverse central $\mathcal C$ est un faisceau en algèbres de Lie semi-simples.}\\





Ce th\'eor\`eme  est en fait un avatar feuillet\'e du  r\'esultat suivant obtenu par Nadel:\\

{\bf Theor\`eme 1.2}  (\cite{nad}) {\it Soit $M$ une vari\'et\'e complexe compacte et connexe dont le fibr\'e canonique est ample. Alors la composante neutre du groupe d'automorphismes $Aut(\tilde M))$ de son rev\^etement universel ${\tilde M}$ est semi-simple.}\\

Tel qu'il est articul\'e, cet \'enonc\'e sous-entend que $Aut(\tilde M))$ est muni d'une structure de groupe de Lie. C'est effectivement vrai et ceci r\'esulte du fait qu'il agit par isom\'etries sur ${\tilde M}$ pour la m\'etrique K\"ahler Einstein relev\'ee. Dans ce contexte, la courbure de Ricci est strictement n\'egative et ceci pr\'ecise les liens avec notre travail. Si $Aut^0(\tilde M))$ est susceptible d'\^etre riche, {\it a contrario}, les automorphismes de $M$ sont en nombre fini, ph\'enom\`enes bien connus lorque $M$ est une surface de Riemann de genre sup\'erieur ou \'egal \`a deux. \\

En fait, ces deux propri\'et\'es, \'a priori antagonistes sont  toutes les deux des cons\'equences des formules de Weitzenb\"ock-Bochner jointes \`a un lemme \'el\'ementaire mais fondamental de Nadel ({\it loc.cit}) concernant les alg\`ebres de Lie ab\'eliennes de champs de Killing et ces techniques s'adaptent parfaitement \`a la situation qui est la n\^otre.\\
Signalons que l'utilisation nous faisons du principe de Bochner est tr\`es similaire \`a celle de Wu et Zheng dans \cite{wu} ( voir aussi \cite{fr}). Chez Nadel ({\it loc.cit}) cet argument est seulement mentionn\'e et ce dernier fait plut\^ot usage de la stabilit\'e du fibr\'e tangent.\\





\vskip 5 pt
On obtient une formulation plus pr\'ecise du th\'eor\`eme 1.1 en le sp\`ecialisant \`a la classe de feuilletages pr\'esent\'ee ci-dessous.\\

\noindent Soit $M$ une vari\'et'e k\"ahlerienne compacte munie d'un feuilletage holomorphe r\'egulier $\mathcal F$. Notons respectivement $T_{\mathcal F}$ et $N_{\mathcal F}$ les fibr\'es tangent et normal de $\mathcal F$ munis de leur structure holomorphe usuelle. Supposons de plus que  l\`a classe canonique de $T\mathcal F$ est num\'eriquement triviale,i.e: 

    $$c_1(T_{\mathcal F})=0$$

Dans \cite{to} est donn\'ee une description assez pr\'ecise de tels feuilletages en codimension 1. En codimension quelconque, nous obtenons ici le 

\vskip 5 pt

{\bf Th\'eor\`eme 1.3} {\it Soit $\mathcal F$ un feuilletage holomorphe r\'egulier comme ci-dessus. On suppose en outre que $c_1(N_{\mathcal F})$ est repr\'esent\'ee par une $(1,1)$ forme ferm\'ee semi-n\'egative et v\'erifie ${c_1(M)}^n\not=0$; alors, \'a rev\^etement fini pr\'es, $\mathcal F$ est une fibration localement triviale au dessus d'une vari\'et\'e $V$ \`a premi\`ere classe de Chern quasi-n\'egative (en particulier, les feuilles de $\mathcal F$ sont ferm\'ees}).\\

\medskip

{\bf Corollaire 1.1} {\it Sous les hypoth\`eses du th\'eor\`eme 1.2, $\mathcal F$
  co\"incide avec la fibration de d'Itaka-Kodaira de $M$.}\\


{\bf Th\'eor\`eme 1.4} {\it Soit $\mathcal F$ un feuilletage holomorphe r\'egulier de codimension $n\leq 2$ \`a classe canonique num\'eriquement triviale sur $M$ k\"ahler compacte. On suppose en outre que $c_1(N_{\mathcal F})$ est repr\'esent\'ee par une $(1,1)$ forme ferm\'ee semi-n\'egative; 
alors il existe  un feuilletage holomorphe $\mathcal G$ de  contenant ${\mathcal F}$ et tel que $c_1(N\mathcal G)=c_1(N\mathcal F)$ et les feuilles de $\mathcal G$ soient fermées.}\\

  Nous conjecturons que cet \'enonc\'e reste valide en laissant tomber l'hypoth\`ese sur la codimension. Comme nous le rappelons en section 4, c'est effectivement le cas quand $\mathcal F$ est d\'efini par l'action localement libre d'un groupe de Lie complexe et ceci r\'esulte simplement des travaux de Lieberman (\cite{li}) relatifs au groupes d'automorphismes des vari\'et\'es k\"ahleriennes.\\

En dehors des ingr\'edients d\'ej\`a mentionn\'es, la d\'emonstration du th\'eor\`eme 1.1 passe par l'obtention d'un lemme $dd^c$ pour les courants basiques. C'est l'objet de la section qui suit.



\bigskip 
\section{ Courants basiques}
\vskip 20 pt

Soit ${\mathcal F}$ un feuilletage de dimension $n$ sur une vari\'et\'e $M$  orientée de dimension (réelle) $N$. 

On notera $\Gamma {\mathcal F}$ le $\mathcal C^\infty (M)$ module des champs de vecteurs sur $M$ tangent \`a ${\mathcal F}$ et on désignera par $\Omega_c^*(M)$ (o\'u simplement $\Omega^*(M)$ quand $M$ est compacte) le complexe des formes différentielles sur $M$ à valeurs complexes et \`a support compact. 

Un courant $T\in {\Omega_{\CC}^*(M)}'$ de degr\'e $r$ est dit {\it basique} si, pour tout $X\in \Gamma {\mathcal F}$, $i_XT=0$ et $L_XT=0$, le produit intérieur $i_X$ et la dérivée de Lie $L_X$ d'un $r$ courant $T$  étant définis par dualité comme suit:

Pour tout champ de vecteurs $X$ sur $M$  et pour toute forme $\alpha\in {\Omega_c^{N-r}(M)}'$,

  $$\langle i_XT,\alpha\rangle={(-1)}^r\langle T,i_X\alpha\rangle$$

 $$\langle L_XT,\alpha\rangle={(-1)}^r\langle T,L_X\alpha\rangle.$$

On notera ${\mathcal C}_{\mathcal F}^*(M)$ le complexe des courants basiques sur $M$.
 
Une $r$  forme $\eta$ définit un courant $T_\eta$ de degr\'e $r$, dit {\it régulier} par la formule habituelle

   $$\langle T_\eta,\alpha\rangle =\int_M \eta\wedge \alpha$$ o\`u $\alpha$ est une $N-r$ forme \'a support compact.

On peut vérifier (\cite{ab}) que $T_\eta$ est basique si et seulement si $\eta$ est basique.\\

Quand $\mathcal F$ est transversalement holomorphe de codimension complexe $r={N-n\over 2}$, un courant basique sera dit par ailleurs de {\it bidegr\'e}  $(p,q)$ s'il appartient au dual topologique de $\Gamma({{N^*\mathcal F}}^{r-p,r-q}\otimes \Lambda^n{T_{\CC}^*\mathcal F})\cap \Omega_c^*(M)$. On note ${\mathcal C}_{\mathcal F}^{p,q}(M)$ l'ensemble de tels courants. On obtient ainsi une d\'ecomposition des courants basiques de degr\'e $r$,

   $${\mathcal C}_{\mathcal F}^{r}(M)=\bigoplus_{p+q=r}\mathcal C_{\mathcal F}^{p,q}(M)$$ et un double complexe

$$ (\mathcal C_{\mathcal F}^{*,*}(M),\partial,\overline{\partial}),$$

\noindent Les op\'erateurs $\partial$ et $\overline{\partial}$ s'\'etendent aux courants basiques comme suit:\\

si $T\in {\mathcal C}_{\mathcal F}^{r}(M)$, on pose pour tout $\omega\in \Omega_c^{N-r}(M)$

$$ \langle\partial T,\omega\rangle={(-1)}^{r+1}\langle T,\partial\omega\rangle$$
et
$$ \langle\overline{\partial} T,\omega\rangle={(-1)}^{r+1}\langle T,\overline{\partial}\omega\rangle$$

Le cadre consid\'er\'e \`a partir de maintenant est celui d'une vari\'et\'e {\it compacte} $M$ munie d'un feuilletage {\it k\"ahl\'erien} de dimension $p$ et de codimension complexe $n$. nous supposerons en outre que $\mathcal F$ est {\it homologiquement orientable}. Cette derni\`ere condition implique, suivant un r\'esultat de Xos\'e Masa (\cite{ma}), que $\mathcal F$ est {\it minimalisable}. Ceci revient  \`a dire qu'il existe sur $M$ une m\'etrique riemannienne $\tilde g$ dont le   volume sur les feuilles est la restriction d'une  forme $\chi\in \Omega^p(M)$  relativement fermée:\\

\noindent pour tous champs de vecteurs $X_1,...,X_{p+1}$ tels que $X_i\in\Gamma (T\mathcal F), i=1,...,p$, on a

$$ d\chi(X_1,...,X_{p+1})=0.$$

On choisira par la suite $\tilde g$ de telle sorte que la m\'etrique induite sur $\nu\mathcal F$ soit pr\'ecis\'ement $g$.\\
L'opérateur $\overline{*}$ est alors reli\'e \'a l'op\'erateur $*$ ordinaire (associ\'e \`a $\tilde g$) par la formule 

$$ *u= \overline{*}u\wedge \chi\leqno (2.1)$$

Pour toute forme $\alpha,\beta\in \Omega^r(M/\mathcal F)$, on obtient donc que 

    $$(\alpha,\beta)=\int_M\alpha\wedge\overline{*}\beta\wedge\chi$$
o\`u pour tout $\xi,\eta\in\Omega^s(M)$, $(\xi,\eta)=\int_M\xi\wedge *\eta$ est le produit scalaire hermitien induit par $\tilde g$ sur $\Omega^*(M)$.\\

Compte-tenu de la relation $(2.1)$, on a par ailleurs que pour tout $u\in \Omega^r(M/\mathcal F)$ et pour tout $\xi\in\Omega^{N-r} (M)$, 

$$\langle T_{\overline{*}u},\xi\rangle=\int_M \overline{*}u\wedge\xi=\pm\int_M*(u\wedge\chi)\wedge\xi=\pm\overline{\langle T_u,{*\xi}\wedge\chi\rangle},$$
\noindent ce qui permet d'expliciter l'extension de $\overline{*}$ (et par suite de $\Delta$) \'a l'espace des courants basiques.\\
On peut montrer qu'on obtient, relativement au produit scalaire pr\'ec\'edent, une {\it d\'ecomposition de Hodge basique} pour l'op\'erateur auto-adjoint $\Delta$ (voir \cite{elk}, \cite{eg}), i.e., une d\'ecomposition en sous-espace orthogonaux du type

$$\Omega^{p,q}(M/\mathcal F)={\mathcal H}^{p,q}\oplus Im\ \Delta$$
o\`u ${\mathcal H}^{p,q}=$ker$\Delta$ est de dimension {\it finie} (par abus de langage, on note encore $\Delta$ la restriction du laplacien basique \'a $\Omega^{p,q}(M/\mathcal F))$.\\
Soit $\omega\in \Omega^{p,q}(M/\mathcal F)$ et $H(\omega)$ sa projection orthogonale sur ${\mathcal H}^{p,q}$; consid\'erons l'op\'erateur de Green $G:\Omega^{p,q}(M/\mathcal F)\rightarrow\Omega^{p,q}(M/\mathcal F)$ (on rappelle que $G\omega$ est l'unique forme telle que $G\omega\in Im\ \Delta$ et $\Delta G\omega=\omega-H(\omega))$.\\

Comme dans le cas classique, 

$$\omega=H(\omega)+\Delta G =H(\omega)+\partial\partial^{\overline{*}}G\omega +\partial^{\overline{*}}\partial G\omega= H(\omega)+{\overline{\partial}}\ {\overline{\partial}}^{\overline{*}}G\omega +{\overline{\partial}}^{\overline{*}}\ {\overline{\partial}} G\omega $$
et $G$ commute aux diff\'erents op\'erateurs $d,\partial,\overline{\partial},\overline{*}.$\\

  Dans ce qui suit, on veut \'etablir qu'il existe une d\'ecomposition similaire pour $T\in \mathcal C_{\mathcal F}^{p,q}(M)$. Pour ce faire, on va utiliser qu'il existe une suite de formes basiques $(u_n)$ de bidegr\'e $(p,q)$ telle que la suite de courant basiques $(T_{u_n})$ converge vers $T$ (\cite{al}). En effet, soit $(u_n)$ une telle suite; pour tout entier $n$ , on a 

$$T_{u_n}=T_{H(u_n)}+T_{\Delta Gu_n}.$$

Puisque ${\mathcal H}^{p,q}$ est de dimension finie, il est clair que $(T_{H(u_n)})$ converge vers $T_{H(T)}$ o\`u $H(T)\in={\mathcal H}^{p,q}$ est l'unique forme harmonique telle que pour tout $h\in {\mathcal H}^{n-p,n-q}$, on ait

$$\int_MH(T)\wedge h\wedge\chi=\langle T,h\wedge\chi\rangle$$

Par suite, $(T_{\Delta G_{u_n}}=\Delta T_{G_{u_n}})$ converge vers un courant basique $T'$.\\
Nous affirmons que la suite de courants $(T_{G_{u_n}})$ est {\it convergente}, auquel cas $T'=\Delta T''$ avec $T''=lim_{n\rightarrow +\infty} T_{G_{u_n}}\in {\mathcal C}_{\mathcal F}^{p,q}(M)$. Ceci se prouve aisément en invoquant le r\'esultat suivant:\\
\medskip

{\bf Th\'eor\`eme 2.1} (\cite{kt}) {\it Il existe un op\'erateur elliptique $$D:\Omega^*(M)\rightarrow \Omega^*(M)$$ pr\'eservant le degr\'e (non  n\'ecessairement auto-adjoint) tel que $D$ co\"\i ncide avec $\Delta$ en restriction \`a $\Omega^*(M/\mathcal F)$ et tel que la d\'ecomposition orthogonale (donn\'ee par la th\'eorie classique) $$\Omega^*(M)=ker\ D\oplus\ Im\ D^*$$ soit compatible avec celle de $\Omega^* (M/\mathcal F)$; en d'autre termes, $ker\ \Delta\subset ker\ D$ et $Im\ \Delta\subset Im\ D^*$.}\\

Etabissons la convergence de 
$T_{G_{u_n}}$; pour tout $v\in \Omega^{p+q}(M)$, on a 

$$\langle T_{\Delta Gu_n},*v\rangle=\langle T_{Gu_n},*D^*v\rangle.$$
Par ailleurs, pour tout $h\in Ker\ D$, on a

$$\langle T_{Gu_n},*h\rangle=(T_{Gu_n},h)=0,$$

\noindent en cons\'equence de quoi $(T_{Gu_n})$ converge (faiblement) vers le courant $T^{''}$ d\'efini par\\

$\langle T^{''},*D^**v\rangle=lim_{n\rightarrow+\infty}\langle T_{\Delta Gu_n},*v\rangle$, pour tout $v\in \Omega^{p+q}(M)$ et\\

$\langle T^{''},*h\rangle=0$, pour tout $h\in Ker\ D\cap \Omega^{p+q}(M)$.\\

La s\'erie d'observations pr\'ec\'edentes permet alors d'\'enoncer la
\medskip

{\bf Proposition 2.1} {\it Soit $(\mathcal F,g)$ un feuilletage kahlerien sur $M$ compacte. Soit $T$ un courant basique de bidegr\'e $(p,q)$; alors $\Delta T=0$ si et seulement si $T=T_h$ o\`u $h\in{\mathcal H}^{p,q}$ et on h\'erite d'une d\'ecomposition en somme directe

    $$ {\mathcal C}_{\mathcal F}^{p,q}(M)=ker\ \Delta\oplus Im\ \Delta.$$

L'op\'erateur de Green $G: {\mathcal C}_{\mathcal F}^{p,q}(M)\rightarrow  {\mathcal C}_{\mathcal F}^{p,q}(M)$  associ\'e \'a cette d\'ecomposition (d\'efini similairement au pr\'ec\'edent) commute \`a $d,\partial, \overline{\partial},\overline{*}$.\\

De plus, si $T$ est $d$ exact (i.e, il existe $T'\in{\mathcal C}_{\mathcal F}^{p+q-1}(M)$ tel que $dT'=T$), on a alors

$$T=2\partial\partial^{\overline{*}} GT=2\overline{\partial}\ \overline{\partial}^{ \overline{*}} GT.$$}\\

\medskip

{\bf Corollaire 2.1} {\it Soit $T\in{\mathcal C}_{\mathcal F}^{p,q}(M)$ $d$ exact; alors $T$ est $\partial\overline{\partial}$ exact; c'est-\`a-dire qu'il existe un courant $T'\in {\mathcal C}_{\mathcal F}^{p-1,q-1}(M)$ tel que 

$$T=\sqrt{-1}\partial\overline{\partial}T'.$$
De plus, on peut choisir $T'$ r\'eel lorsque $T\in {\mathcal C}_{\mathcal F}^{p,p}(M)$ est r\'eel.}

\section{Preuve du th\'eor\`eme 1.1}
\vskip 20 pt

On rappelle que, par hypoth\`ese,  $({\mathcal F}, g)$ d\'esigne  un feuilletage transversalement kahlerien et homologiquement orientable dont la courbure de Ricci transverse $\rho$ est cohomologue (dans  $\Omega^*(M/\mathcal F)$) \`a une $(1,1)$ forme r\'eelle basique quasi-n\'egative $\kappa$. Quitte \`a se placer sur un double rev\^etement, on peut \'egalement supposer que $\mathcal F$ est orientable.\\
Apr\`es  r\'esolution de l'\'equation de Monge-Amp\`ere basique (\cite{elk}), on se ram\`ene au cas o\`u   $\kappa=\rho$. Un des points cl\'es de la d\'emonstration est certainement le lemme suivant qui est de nature purement locale:\\
\vskip 10 pt
{\bf Lemme 3.1}. (\cite{nad})- {\it Soient $U$ un ouvert connexe de ${\CC}^n$ muni d'une m\'etrique kahlerienne $g$ (non n\'ecessairement compl\`ete) et $X_1,...,X_l$ des champs de vecteurs holomorphes sur $U$ commutant 2 \`a 2 tels que leur partie r\'eelle $Re(X_i),\ i=1,...,l$ soit un champ de Killing. \\
On suppose de plus que les $X_i,i=1,...,n$ forment une famille libre sur ${\RR}$ mais sont m\'eromorphiquement d\'ependants; alors la courbure de Ricci $\rho$ n'a de signe d\'efini en aucun point de $U$.}\\ 

Supposons par l'absurde que $\mathcal C$ ne soit pas un faisceau d'alg\`ebre semi-simples. Chaque germe d'alg\`ebre ${\mathcal C}_m,\ m\in M$ admet donc un radical r\`esoluble $R,m$ {\it non trivial}, lequel est unique. La collection des $R_m$ d\`efinit donc un faisceau localement d'alg\`ebre de Lie r\`esolubles. Similairement, on peut consid\'erer le sous-faisceau localement constant $\mathcal A$ dont la fibre, en chaque point $m$, est la derni\`ere alg\`ebre d\'eriv\'ee non triviale.\\
Par ``complexification'', on h\'erite alors d'un fasceau localement constant ${\mathcal A}_{\CC}$ d'alg\`ebres de Lie de champs de vecteurs holomorphes. Plus pr\'ecis\'ement, ${\mathcal A}_{\CC}$ est engendr\'ee sur ${\CC}$ \`a partir des sections locales $X-iJX$ o\`u $X$ est une section locale de ${\mathcal A}$ et $J$ la structure complexe transverse. Remarquons que $d:=Dim_{\CC}{\mathcal A}_{\CC}\leq Dim_{\RR}{\mathcal A}$.\\

Par \'evaluation des champs de vecteurs, on r\'ecup\`ere un morphisme de faisceaux:

$$ ev:{\mathcal A}_{\CC}\rightarrow N_{\mathcal F}^{1,0}$$

qui induit un morphisme:

  $$\bigwedge^d ev:det {\mathcal A}_{\CC}=\bigwedge^d  {\mathcal A}_{\CC}\rightarrow        \bigwedge^d  N_{\mathcal F}^{1,0}$$.

Compte-tenu du lemme 3.1  et de l'hypoth\`ese $\rho$ quasi-n\'egatif'', ce morphisme n'est pas trivial et on r\'ecup\`ere par ce proc\'ed\'e une section {\it holomorphe non nulle} $s$ du fibr\'e $E:=  \bigwedge^d  N_{\mathcal F}^{1,0}\otimes {(det {\mathcal A}_{\CC})}^{-1}$ .\\

 On peut donc trouver un recouvrement ouvert $\mathcal U=\{U_i\}$ de $M$ de telle sorte qu'en restriction \'a chaque $U_i$, $s$ soit repr\'esent\'e par une section holomorphe $s_i$ de $\bigwedge^d  N_{\mathcal F}^{1,0}$ telle que sur   
 chaque intersection $U_i\cap U_j$, 

 $$s_i=g_{ij}s_j$$
o\`u $(g_{ij})$ est un cocycle multiplicatif {\it localement constant}.

Soit $g_d$ la m\'etrique induite par $g$ sur le fibr\'e  $\bigwedge^d  N_{\mathcal F}^{1,0}$. Sur $ U_i\cap U_j$, on a

  $$g_d(s_i,\overline{s_i})=|g_{ij}|^2 g_d(s_j,\overline{s_j}).$$
Par ailleurs,il est clair que pour tout indice $i$,  $\varphi_i=\log (g_d(s_i,\overline{s_i})\in L_{loc}^1 (U_i)$ et ces fonction sont localement constantes sur les feuilles de $\mathcal F_{|U_i}$.
Par recollement, la collection des  $\sqrt{-1}$ $\partial\overline {\partial}\varphi_i=$ (bien d\'efinis au sens de courants) produit un courant basique $T$ de bidegr\'e $(1,1)$ $d$ exact (en cohomologie basique). Il r\'esulte alors du corollaire 2.1 que 
$$ T=\sqrt{-1}\partial\overline {\partial} u$$ o\`u $u\in {\mathcal C}_\mathcal F^{0}(M)$ est une distribution basique {\it r\'eelle}. Sur tout ouvert $U_i$, on obtient ainsi que 

       $$u-T_{\varphi_i}=T_{Log{|H_i|}^2}$$ o\`u $H_i$ est une fonction holomorphe basique (pour le feuilletage restreint) ne s'annulant pas sur $U_i$. Il en r\'esulte que pour tout $i$,

   $$H_is_i=h_{ij}H_js_j$$
o\`u les $h_{ij}$ sont localement constants de {\it module 1}.\\

Compte-tenu de ce qui pr\'ec\`ede, les ${||H_is_i||}^2$ se recollent sur les intersections en une fonction positive $\xi$ telle que, suivant la formule de Weintzenb\"ock, on ait

$$\bigtriangleup\xi\gneqq, 0$$

\noindent ce qui contredit le th\'eor\`eme de Hopf.\qed




\medskip

{\bf Remarque 3.1:} Les m\^emes d\'emonstrations et la conclusion du th\'eor\`eme 1.1 conclusions restent valides en remplaçant $\mathcal C$ par n'importe quel faisceau localement constant de champs de killing transverses qui sont des parties r\'eelles de champs basiques holomorphes locaux.

\bigskip
\section{Feuilletages \`a classe canonique num\'eriquement triviale}
\vskip 20 pt

Dans cette section, on d\'esignera par ${\mathcal F}$ un feuilletage holomorphe r\'egulier sur une vari\'et\'e $(M,g)$ k\"ahleriene compacte de dimension complexe $m$ qui satisfait une partie des hypoth\'eses du th\'eor\`eme 1.1, \'a savoir $\mathcal F$ est \`a classe canonique num\'eriquement triviale, i.e, $c_1(T_{\mathcal F})=0$; on supposera en revanche la classe anticanonique $c_1(M)$  repr\'esent\'ee par une $(1,1)$ forme ferm\'ee $\eta$ seulement {\it semi-n\'egative} (et \`a priori non basique).\\

{\bf 4.1  Preuve du th\'eor\`eme 1.3}\\

Elle r\'esulte de la s\'erie d'observations suivantes.\\

{\bf Lemme 4.1} {\it Il existe sur $M$ une m\'etrique kahlerienne pour laquelle le feuilletage ${\mathcal F}$ (plus exactement le fibr\'e $T_{\mathcal F}$) est parall\`ele. En particulier, $\eta$ est une forme basique de $\mathcal F$.}\\
\medskip

{\it Preuve}-. Suivant le fameux th\'eor\`eme de Yau (\cite{ya}),  on peut munir $M$ d'une m\'etrique k\"ahlerienne  dont la courbure de Ricci est pr\'ecis\'ement $\eta$. Par hypoth\`ese, il existe  une section holomorphe {\it non triviale} du fibr\'e en droites $\bigwedge^{m}T_{\mathcal F}\otimes E$ o\`u $E$ est un certain fibr\'e {\it plat}. Le lemme r\'esulte alors du principe de Bochner tel qu'il est par exemple utlilis\'e dans la section pr\'ec\'edente (voir aussi \cite{to}).\qed\\
\medskip

Notons dor\'enavant par $g$ la m\'etrique ainsi produite. On h\'erite par parall\'elisme d'une d\'ecomposition {\it orthogonale} et {\it holomorphe} du fibr\'e tangent de $M$:

  $$TM=T_{\mathcal F}\oplus {T_{\mathcal F}}^\perp$$
o\`u  ${T_{\mathcal F}}^\perp=T_{\mathcal G}$ est le fibr\'e tangent d'un feuilletage holomorphe $\mathcal G$ \'egalement parall\`ele.\\
Ce scindage, de nature infinit\'esimale, donne lieu par le th\'eor\`eme de de Rham \'a une d\'ecomposition  du rev\^etement universel $\tilde M$ en produit  kahlerien  pour la m\'etrique relev\'ee $h=h_1\oplus h_2$:

       $$\tilde M=(V,h_1)\times (V^\perp,h_2).$$
Par construction, la m\'etrique $h_1$  est Ricci plate, ce qui va permettre de pr\'eciser la structure de $V$.\\
\medskip

{\bf Lemme 4.2} {\it La variété ${V}$ se scinde (holomorphiquement et
isométriquement) sous la forme ${\CC}^k\times N$ où $N$ est une
variété de Calabi-Yau.}\\
\medskip

La d\'emonstration de ce lemme suit celle qui est pr\'esent\'ee dans \cite{to}; nous la reproduisons par commodit\'e pour le lecteur.

Elle repose notamment sur le résultat suivant de Cheeger et Gromoll.
\vskip 10 pt {\bf Théorème 4.1} (\cite{ch})
\vskip 5 pt
{\it Soit $M$ une variété lisse simplement connexe admettant une
métrique riemannienne complète de courbure de Ricci positive ou
nulle; alors $M$ se décompose isométriquement sous la forme ${\RR}^k\times N$ où $N$ est une variété ne contenant pas de droite
géodésique.} \vskip 5 pt Rappelons qu'une droite géodésique est
une géodésique

$$ \gamma:]-\infty, +\infty[\rightarrow M$$
telle qu'à chaque instant $t,t'\in{\RR}$, le segment $\gamma_{|[t,t']}$ soit
minimal. \medskip

  {\it Preuve du lemme 4.2}-. Puisque la d\'ecomposition de De Rham d'une variété k\"ahlerienne
co\"\i ncide avec celle de la variété réelle sous-jascente,  on obtient que ${V}$ se scinde
holomorphiquement et isométriquement sous la forme ${\CC}^l\times N$.

Par construction, le groupe fondamental de $M$ agit diagonalement
par isom\'etries sur le produit ${\tilde M}={\CC}^l\times
N\times V^\perp$. Il reste à prouver que $N$ est compacte;
supposons par l'absurde que ce ne soit pas le cas. La
contradiction recherchée résulte alors de l'argument suivant
repris {\it verbatim} à Cheeger et Gromoll ({\it loc.cit}). Soit
$K$ un domaine fondamental pour l'action de $\pi_1(M)$ sur
${\tilde M}$. C'est un compact et sa projection $\rho (K)$ sur $N$
est donc un compact dont l'orbite par $\rho (\pi_1(M))$ est $N$
toute entière. Puisque $N$ est suppos\'ee non compacte, il existe en
un point $p\in \rho (K)$ un rayon géodésique
$\gamma:[0,\infty[\rightarrow V$(i.e une demi-droite g\'eod\'esique)
tel que $\gamma (0)=p\in \rho (K)$. Soit $g_n$ une suite
d'isom\'etries de $\rho (\pi_1(M))$ telle que $g_n(\gamma(n))=p_n\in
\rho (K)$. La g\'eod\'esique $\gamma_n$ de $N$ d\'efinie par
$\gamma_n(0)=p_n$ et ${\gamma_n}^{'}(0)=dg_n({\gamma}^{'}(n))$ est
donc un rayon g\'eod\'esique en restriction à $[-n, \infty[$. Par
compacit\'e, on peut
 extraire une sous-suite $n_i$ de sorte que $p_{n_i}$ et
${\gamma_{n_i}}^{'}(0)$ convergent respectivement vers $p_0\in
\rho (K)$ et $v_0\in T_{p_0}^1 (N)$. Par suite, la
g\'eod\'esique passant en $p_0$ à la vitesse $v$ est une droite
g\'eod\'esique.\qed

\vskip 5 pt
{\bf Lemme 4.3} {\it Supposons de plus que les feuilles de $\mathcal F$ sont ferm\'ees, alors il existe un rev\^etement holomorphe fini de $M$ tel que le feuilletage pull-back soit une fibration holomorphe localement triviale.}
\medskip

{\it Preuve du lemme 4.3-.} 
 Le  groupe fondamental $\pi_1(M)$ agit sur $\tilde M={\CC}^k\times N\times V^\perp.$ diagonalement par isom\'etries biholomorphes. Quitte \`a substituer \`a $M$ un rev\^etement fini, on peut supposer que $\pi_1(M)$ agit trivialement sur le facteur Calabi-Yau $N$ (car son groupe d'isom\'etries est fini (\cite{be})). Le lemme 4.3 sera \'etabli si l'on arrive \`a extraire de $\pi_1(M)$ un sous groupe d'indice fini dont l'image par la projection $\rho:\pi_1(M)\rightarrow V\perp $ est sans torsion. Cela s'effectue en suivant pas \`a pas le raisonnement men\'e dans \cite{wu}, p.279.\qed.
\medskip

    \medskip        

D'apr\'es ce qui pr\'ec\`ede, le th\'eor\`eme 1.3 sera d\'emontr\'e une fois acquis le \\

{\bf Lemme 4.4} {\it Supposons que $c_1(M)$ soit repr\'esent\'ee par une forme quasi-n\'egative; alors l'action de $\pi_1(M)$ sur ${V^\perp}$ est
discrète.} \vskip 10 pt
 Puisque le groupe des isom\'etries de $N$ est fini (\cite{be}), on est ramen\'e à consid\'erer un groupe $H$ dont l'action est  libre,
 discr\`ete, diagonale, cocompacte et isom\'etrique sur le produit $ {\CC}^l\times {V^\perp}$. Il reste à voir que cette action reste
  {\it discrète} sur le second facteur, auquel cas chaque feuille de $\mathcal F$ sera {\it fermée}, ce qui ach\`evera la preuve du th\'eor\'eme 1.3.\\

\medskip

{\it Preuve du lemme 4.4}-. Puisque l'action du groupe fondamental est diagonale, on r\'ecup\`ere un  morphisme  $\rho$ de $\pi_1(M)$ vers le groupe de Lie $Isom  ({V^\perp})$  des isom\'etries {\it holomorphes} de ${V^\perp}$. Soit $G$ la composante neutre de l'adh\'erence de $\rho (\pi_1(M)$ dans $Isom  ({V^\perp})$ et $\mathfrak g$ son alg\`ebre de lie. Suivant un r\'esultat d\^u \`a Eberlein, (\cite{eb}), $\mathfrak g$ est r\'esoluble; elle est par ailleurs semi-simple compte tenu du th\'eor\`eme 1.1\footnote{On peut le d\'eduire plus simplement du principe de Bochner usuel (cf.\cite{wu}) mais il nous a sembl\'e interessant d'\'etablir un r\'esultat plus g\'en\'eral}.Par suite,  $\mathfrak g=\{0\}$ et le lemme est d\'emontr\'e.\qed\\








{\it Preuve du corollaire 1.1-. }
 On peut supposer (lemme 4.3) que $\mathcal F$  est une fibration localement triviale. D'après \cite{SS}, la base $V$ de cette fibration est de type g\'en\'eral D'autre part, les fibres sont \`a premi\`ere classe Chern nulle  et ceci implique qu'il exite un entier $n$ non nul tel que $K_{\mathcal F}^{\otimes n}$ soit trivial. Le corollaire r\'esulte alors simplement du fait que $K_{M'}=\pi^*K_V\otimes K_{\mathcal F}$. \qed \\

{\bf 4.2  Preuve du th\'eor\`eme 1.4}\\

Examinons d'abord le cas o\`u le feuilletage en question est donn\'ee par une action localement libre d'un groupe de Lie complexe $G$; le th\'eor\`eme 1.4 est une cons\'equence directe, sans hypoth\`ese sur la codimension, d'un r\'esultat d\'emontr\'e id\'ependamment par Lieberman et Fujiki. Ce r\'esultat stipule que la composante connexe $Aut_0(M)$ du groupe des automorphismes d'une vari\'et\'e k\"ahlerienne compacte $M$ jouit des m\^emes propri\'et\'es qu'un groupe alg\'ebrique; en particulier, le feuilletage holomorphe $\mathcal G$ (peut- \^etre singulier) d \'efini par les orbites de  $Aut_0(M)$ (o\`u celles de l'adh\'erence de Zariski de $G$)  est \`a feuilles ferm\'ees (dans le compl\'ementaire du lieu singulier \'eventuel). En courbure de Ricci semi-n\'egative, tout champ de vecteur holomorphe est de plus parall\`ele d'apr\`es le principe de Bochner. Ceci nous assure que $\mathcal G$ est d\'efini par une action localement libre et qu'en particulier $c_1(T\mathcal G)=0$.\\

Le  th\'eor\`eme 1.4 a d\'ej\`a \'et\'e \'etabli lorsque $c_1(M)$ est de rang maximal $n=$codim$\mathcal F$ quelque part ($\mathcal G=\mathcal F$ convient) et le seul cas consistant \`a traiter est \'evidemment $n=2$.

Soit $G$ la composante neutre de $\overline{\rho (\pi_1(M)}$ dont on rappelle que  l'alg\`ebre de Lie $\mathfrak g$ est  r\'esoluble. Les conclusions du th\'eor\`eme sont bien sur v\'erifi\'ees si celle-ci est triviale. On supposera donc que ce n'est pas le cas.

D\'esignons par $\mathfrak h$ sa derni\`ere  alg\`ebre de Lie d\'eriv\'ee non triviale; elle est globalement $\pi_1(M$ invariante et admet pour base $p$ champs de Killing $v_1,...,v_p$, chaque $v_i$ \'etant la partie r\'eelle d'un champ de vecteur {\it holomorphe} $X_i$ sur $V^{\perp}$. Par d\'efinition de $\mathfrak h$, les champs $V_i$, et par suite les champs $X_i$ commutent deux \`a deux. Notons $p'$ la dimension complexe de l'alg\`ebre de Lie ${\mathfrak h}_{\CC}$ engendr\'ee sur ${\CC}$ par les $X_i$. Il est clair que $p'\leq p$. 

Supposons d'abord que $X_1\wedge X_2\not=0$. On en d\'eduit
ais\'ement, par ab\'elianit\'e, que $\{ X_1, X_2\}$ forme une base de ${\mathfrak h}_{\CC}$.

On a dautre part une repr\'esentation 

$$\pi_1(M)\rightarrow GL({\mathfrak h}_{\CC})$$

\noindent qui donne donc lieu \`a un fibr\'e holomorphe {\it plat} de rang 2 sur la vari\'et\'e $M$. On r\'ecup\`ere par ce proc\'ed\'e une section non triviale $s$ du fibr\'e
$\bigwedge^n TM\otimes{\bigwedge^2 E}^{-1}$. Par les techniques de Bochner d\'ej\`a mentionn\'ees, $s$ est parall\`ele et ne s'annule donc pas; en cons\'equence de quoi on a $c_1(M)=0$ et le th\'eor\`eme 1.4 est montr\'e ( $\mathcal G$ est ici de codimension nulle!).

Similairement, on peut construire sous l'hypoth\`ese $p'=1$  un feuilletage $\mathcal G$ de codimension 1 satisfaisant les conclusions du th\'eor\`eme 1.4..

Il reste donc \`a consid\'erer le cas $p\geq2$ avec $X_i\wedge X_j=0$ pour tout couple d'indices $(i,j)$. Quitte \`a renum\'eroter les indices, on peut supposer que $X_1, ....,X_{p'}$ forment une base de ${\mathfrak h}_{\CC}$.

Soit $m\in V^\perp$ tel que $X_1(m)\not=0$ et $(z_1,z_2)$ un syst\`eme de coordonnées holomorphes locales en $m$

$$ X_1={\partial\over\partial z_1}$$
On a alors

 $$X_i=a_i(z_2),{\partial\over\partial z_1}i=2,...,p',$$
$a_i$ holomorphes non constantes.

En exploitant le fait que les $v_k, k=1,...,p'$ sont de Killing pour la m\'etrique k\"alherienne 

           $$g^\perp= \sum_{i,j=1,2}g_{ij}dz_i\otimes d\overline{z_j},$$

\noindent on obtient facilement que $(a_k-\overline{a_k})g_{11}$ est constant pour tout $k$. Il en r\'esulte que $p'=2$ et que le champ de bivecteur $Re\ X_1\wedge Re\ X_2$ est de norme constante non nulle. Par suite, $Re\ X_1$ et $Re\ X_2$ sont les g\'en\'erateurs d'un feuilletage holomorphe $\tilde{\mathcal G}$ sur $V^\perp$ d\'efini localement par l'\'equation $\{dz_2=0\}$. Notons que $\tilde{\mathcal G}$ est hermitien, car d\'etermin\'e par l'action de champs de Killing (la m\'etrique transverse invariante est celle induite sur ${T\tilde{\mathcal G}}^\perp$ par $h_2$). Ce feuilletage est invariant sous l'action de $\overline{\rho (\pi_1(M)}$, vu que $\mathfrak h$ l'est. De ces quelques remarques, il r\'esulte que le feuilletage $\rho^*(\tilde{\mathcal G})$ se redescend sur $M$ en un feuilletage $\mathcal G$ holomorphe de codimension 1 complexe hermitien (et donc k\"ahlerien).
\medskip

{\bf Lemme 4.5} {\it Les feuilles de $\mathcal G$ sont ferm\'ees}\\

{\it Preuve-.} Supposons par l'absurde que ce n'est pas le cas. Soit $\{X_1,X_2=fX_1\}$ une base  de $\mathfrak h_{\CC}$ D'apr\'es ce qui pr\'ec\'ede, on a $X_2==fX_1$ o\'u $f$ est une fonction holomorphe dont la partie imaginaire ne s'annule pas et telle que $X_1(f)=0$. Pour tout $g\in\rho(\pi_1(M))$, il existe manifestement une matrice 
$$\begin{pmatrix}a&b \\
           c&d\end{pmatrix}\in Gl_2(\mathbb R)$$
telle que 
$$\left\{\begin{matrix}g_*X_1&=&aX_1+bX_2\\
                 g_*X_2&=&cX_1+dX_2\end{matrix}\right.$$

Par suite, la $(1,1)$ forme $\sqrt{-1}/2\pi {df\wedge d\overline{f}\over {|f-\overline{f}|}^2}$ se redescend sur $M$ en une m\'etrique (\'eventuellement d\'eg\'en\'er\'ee) transverse \'a $\mathcal G$ et invariante par holonomie d\'efinie par une $(1,1)$ forme $\eta$. Le lieu des z\'eros de $\eta$ est une union finie d'hypersurfaces $\bigcup_{i=1}^lH_i$ invariantes par $\mathcal G$. Cette forme $\eta$ permet de munir le fibr\'e normal $N_\mathcal G$ d'une m\'etrique (\'eventuellement singuli\'ere) dont la forme de courbure $\xi$ est telle que 

$$\xi =-\sum_{i=1}^ln_i\delta_{H_i}-\eta$$ o\'u les $n_i$ sont des entiers positifs ou nuls et le symbole $\delta_{H_i}$ d\'esigne le courant d'int\'gration le long de $H_i$. La classe $c_1(N^*\mathcal G)$ est donc pseudo-effective non triviale (i.e repr\'esent\'ee par un courant positif ferm\'e non r\'eduit \'a z\'ero) et si l'on reprend l'analyse men\'ee dans \cite{to}, ceci montre que $\mathcal G$ est en fait minimal; en particulier, $n_i=0$ pour tout $i=1,...,l$. Par suite, $\eta$ est une m\'etrique {\it lisse } \'a courbure constante {\it n\'egative}.  
 Soit ${\mathfrak h}'\supseteq{\mathfrak h}$ l'id\'eal de $\mathfrak g$ constitu\'e des champs tangents \`a $\tilde{\mathcal G}$; l'alg\'ebre de lie r\'esoluble $\mathfrak g/{\mathfrak h}'$  se redescend via $\pi:\tilde M\rightarrow M$ en un faisceau localement constant d'alg\'ebre de lie de champs de killing transverses \'a $\mathcal G$. Ce faisceau est en fait r\'eduit \'a z\'ero  d'apr\'es  la remarque 2.3 et les feuilles de $\mathcal G$ sont donc bien ferm\'ees.\qed\\

Ainsi, $\mathcal G$ est une fibration holomorphe au-dessus d'une surface de Riemann; qui plus est, la description de l'alg\'ebre de Lie ${\mathfrak h}_{\CC}$ nous permet d'affirmer que pour toute feuille $\mathcal L$, on a

$$c_1({K_\mathcal H}_{|\mathcal L})=0.$$

Il existe par suite un entier naturel $l$ tel que ${K_\mathcal G}^{\otimes l}$ soit holomorphiquement trivial; il est en effet facile d'en extraire une section $s$ v\'erifiant par exemple pour toute feuille $\mathcal L$
     $${\sqrt {-1}}^{{(n-1)}^2}\int_\mathcal L {(s\wedge\overline{s})}^{1\over l}= \int_\mathcal L \theta^{n-1}$$

o\'u $\theta$ est une forme de kahler fix\'ee sur $M$.

\end{document}